\documentclass[11pt]{amsart}
\usepackage{amsfonts}

\usepackage{graphicx}
\usepackage{amscd}

\newtheorem{theorem}{Theorem}
\theoremstyle{plain}

\newtheorem{example}{Example}

\newtheorem{remark}{Remark}

\numberwithin{equation}{section}

\input{tcilatex}

\begin{document}
\title[Cuntz-Krieger algebras arising from real quadratic maps]{K-theory for Cuntz-Krieger algebras arising from real quadratic maps}
\author{Nuno Martins}
\address{Departamento de Matem\'{a}tica, Instituto Superior T\'{e}cnico}
\email{nmartins@math.ist.utl.pt}
\author{Ricardo Severino}
\address{Departamento de Matem\'{a}tica, Universidade do Minho}
\email{ricardo@math.uminho.pt}
\author{J. Sousa Ramos}
\address{Departamento de Matem\'{a}tica, Instituto Superior T\'{e}cnico}
\email{sramos@math.ist.utl.pt}
\date{September 19, 2002}

\begin{abstract}
We compute the $K$-groups for the Cuntz-Krieger algebras $\mathcal{O}_{A_{%
\mathcal{K}(f_{\mu })}}$, where $A_{\mathcal{K}(f_{\mu })}$ is the Markov
transition matrix arising from the \textit{kneading sequence }$\mathcal{K}%
(f_{\mu })$ of the one-parameter family of real quadratic maps $f_{\mu }$.
\end{abstract}

\maketitle

\baselineskip0.6cm

Consider the one-parameter family of real quadratic maps $f_{\mu
}:[0,1]\rightarrow \lbrack 0,1]$ defined by $f_{\mu }(x)=\mu x(1-x)$, with $%
\mu \in \lbrack 0,4]$. Using Milnor-Thurston's kneading theory [\textbf{14}%
], J. Guckenheimer [\textbf{5}] has classified up to topological conjugacy,
a certain class of maps which includes the quadratic family. The idea of
kneading theory is to encode information about the orbits of a map in terms
of infinite sequences of symbols and to exploit the natural order of the
interval to establish topological properties of the map. In the following, $I
$ will denote the unit interval $[0,1]$ and $c$ the unique turning point of $%
f_{\mu }$. For $x\in I$, let 
\begin{equation*}
\varepsilon _{n}(x)=\left\{ 
\begin{array}{lll}
-1\text{ } & \text{if } & f_{\mu }^{n}(x)>c\text{,} \\ 
\text{ \ }0 & \text{if} & f_{\mu }^{n}(x)=c\text{,} \\ 
+1 & \text{if} & f_{\mu }^{n}(x)<c\text{.}
\end{array}
\right. 
\end{equation*}
The sequence $\varepsilon (x)=(\varepsilon _{n}(x))_{n=0}^{\infty }$ is
called the itinerary of $x$. The itinerary of $f_{\mu }(c)$ is called the 
\textit{kneading sequence} of $f_{\mu }$ and will be denoted by $\mathcal{K}%
(f_{\mu })$. Observe that $\varepsilon _{n}(f_{\mu }(x))=\varepsilon
_{n+1}(x)$, i.e. $\varepsilon (f_{\mu }(x))=\sigma \varepsilon (x)$ where $%
\sigma $ is the shift map. Let $\sum =\{-1,0,+1\}$ be the alphabet set. The
sequences on $\sum^{\mathbb{N}}$ are ordered lexicographi\-cally. However,
this ordering is not reflected by the mapping $x\rightarrow \varepsilon (x)$%
, because the map $f_{\mu }$ reverses orientation on $[c,1]$. To take this
into account, for a sequence $\varepsilon =(\varepsilon _{n})_{n=0}^{\infty }
$ of the symbols $-1,0,+1$, another sequence $\theta =(\theta
_{n})_{n=0}^{\infty }$ is defined by $\theta _{n}=\prod_{i=0}^{n}$ $%
\varepsilon _{i}$. If $\varepsilon =\varepsilon (x)$ is the itinerary of a
point $x\in I$ then $\theta =\theta (x)$ is called the \textit{invariant
coordinate} of $x$. The fundamental \ observation of J. Milnor and W.
Thurston [\textbf{14}] is the monotonicity of the invariant coordinates: 
\begin{equation*}
x<y\Rightarrow \theta (x)\leq \theta (y)\text{.}
\end{equation*}

Let us now consider only those kneading sequences that are periodic, i.e. 
\begin{equation*}
\begin{array}{ll}
\mathcal{K}(f_{\mu }) & =\varepsilon _{0}(f_{\mu }(c))...\varepsilon
_{n-1}(f_{\mu }(c))\varepsilon _{0}(f_{\mu }(c))...\varepsilon _{n-1}(f_{\mu
}(c))... \\ 
& =(\varepsilon _{0}(f_{\mu }(c))...\varepsilon _{n-1}(f_{\mu }(c)))^{\infty
}\equiv (\varepsilon _{1}(c)...\varepsilon _{n}(c))^{\infty }
\end{array}
\end{equation*}
for some $n\in \mathbb{N}$. The sequences $\sigma ^{i}(\mathcal{K}(f_{\mu
}))=\varepsilon _{i+1}(c)\varepsilon _{i+2}(c)...,$ $i=0,1,2,...$, will then
determine a Markov partition of $I$ into $n-1$ line intervals \{$%
I_{1},I_{2},...,\allowbreak I_{n-1}$\} [\textbf{15}], whose definitions will
be given in the proof of Theorem 1. Thus, we will have a Markov transition
matrix $A_{\mathcal{K}(f_{\mu })}$ defined by 
\begin{equation*}
A_{\mathcal{K}(f_{\mu })}:=(a_{ij})\text{ with }a_{ij}=\left\{ 
\begin{array}{l}
1\text{ if }f_{\mu }(\text{int }I_{i})\supseteq \text{int }I_{j} \\ 
\\ 
0\text{ otherwise.}
\end{array}
\right.
\end{equation*}
It is easy to see that this matrix $A_{\mathcal{K}(f_{\mu })}$ is not a
permutation matrix and no row or column of $A_{\mathcal{K}(f_{\mu })}$ is
zero. Thus, for each one of these matrices and following the work of J.
Cuntz and W. Krieger [\textbf{2}], one can construct the Cuntz-Krieger
algebra $\mathcal{O}_{A_{\mathcal{K}(f_{\mu })}}$. In [\textbf{3}], J. Cuntz
proved that 
\begin{equation*}
K_{0}(\mathcal{O}_{A})\cong \mathbb{Z}^{r}/(1-A^{T})\mathbb{Z}^{r}\text{ \ \
and\ \ \ }K_{1}(\mathcal{O}_{A})\cong \text{ }\ker (I-A^{t}:\mathbb{Z}%
^{r}\rightarrow \mathbb{Z}^{r})\text{.}
\end{equation*}
for a $r\times r$ matrix $A$ that satisfies a certain condition (I) (see [%
\textbf{2}]), which is readily verified by the matrices $A_{\mathcal{K}%
(f_{\mu })}$. In [\textbf{1}] R. Bowen and J. Franks introduced the group $BF(A):=$
$\mathbb{Z}^{r}/(1-A)\mathbb{Z}^{r}$ as an invariant for flow equivalence of
topological Markov subshifts determined by $A$.

We can now state and prove the following.

\begin{theorem}
Let $\mathcal{K}(f_{\mu })=(\varepsilon _{1}(c)\varepsilon
_{2}(c)...\varepsilon _{n}(c))^{\infty }$, for some $n\in \mathbb{N}%
\backslash \left\{ 1\right\} $. Thus, we have 
\begin{equation*}
K_{0}(\mathcal{O}_{A_{\mathcal{K}(f_{\mu })}})\cong \mathbb{Z}_{a}\text{ \ \
with}\ \ \ a=\left|
1+\sum\limits_{l=1}^{n-1}\prod\limits_{i=1}^{l}\varepsilon _{i}(c)\right|
\end{equation*}
and 
\begin{equation*}
K_{1}(\mathcal{O}_{A_{\mathcal{K}(f_{\mu })}})\cong \left\{ 
\begin{array}{ccc}
\{0\} &  & \text{if \ \ }a\neq 0 \\ 
&  &  \\ 
\mathbb{Z} &  & \text{if \ \ }a=0\text{.}
\end{array}
\right.
\end{equation*}

\begin{proof}
Set $z_{i}=\varepsilon _{i}(c)\varepsilon _{i+1}(c)...$ for $i=1,2...$. Let $%
z_{i}^{\prime }=f_{\mu }^{i}(c)$ be the point on the unit interval $[0,1]$
represented by the sequence $z_{i}$ for $i=1,2...$. We have $\sigma (z_{i})=$
$z_{i+1}$ for $i=1,...,n-1$ and $\sigma (z_{n})=z_{1}$. Denote by $\omega $
the $n\times n$ matrix representing the shift map $\sigma $. Let $C_{0}$ be
the vector space spanned by the formal basis $\left\{ z_{1}^{\prime },\ldots
,z_{n}^{\prime }\right\} $. Now, let $\rho $ be the permutation of the set $%
\left\{ 1,\ldots ,n\right\} $, which allows us to order the points $%
z_{1}^{\prime },\ldots ,z_{n}^{\prime }$ on the unit interval $[0,1]$, i.e. 
\begin{equation*}
0<z_{\rho (1)}^{\prime }<z_{\rho (2)}^{\prime }<\cdots <z_{\rho (n)}^{\prime
}<1\text{.}
\end{equation*}
Set $x_{i}:=z_{\rho (i)}^{\prime }$ with $i=1,\ldots ,n$ and let $\pi $
denote the permutation matrix which takes the formal basis $\left\{
z_{1}^{\prime },\ldots ,z_{n}^{\prime }\right\} $ to the formal basis $%
\left\{ x_{1},\ldots ,x_{n}\right\} $. We will denote by $C_{1}$ the $n-1$
dimensional vector space spanned by the formal basis $\left\{
x_{i+1}-x_{i}:i=1,\ldots ,n-1\right\} $. Set 
\begin{equation*}
I_{i}:=[x_{i},x_{i+1}]\text{, \ \ for \ \ }i=1,\ldots ,n-1\text{.}
\end{equation*}
Thus, we can define the Markov transition matrix $A_{\mathcal{K}(f_{\mu })}$
as above. Let $\varphi $ denote the incidence matrix that takes the formal
basis $\left\{ x_{1},\ldots ,x_{n}\right\} $ of $C_{0}$ to the formal basis $%
\left\{ x_{2}-x_{1},\ldots ,x_{n}-x_{n-1}\right\} $ of $C_{1}$. Put $\eta
:=\varphi \pi $. As in [\textbf{7}] and [\textbf{8}], we obtain an
endomorphism $\alpha $ of $C_{1}$, that makes the following diagram
commutative. 
\begin{equation*}
\begin{array}{ccccc}
&  & \eta &  &  \\ 
& C_{0} & \longrightarrow & C_{1} &  \\ 
\omega & \downarrow &  & \downarrow & \alpha \\ 
& C_{0} & \longrightarrow & C_{1} &  \\ 
&  & \eta &  & 
\end{array}
\end{equation*}
We have $\alpha =\eta \omega \eta ^{T}(\eta \eta ^{T})^{-1}$. Remark that if
we neglect the negative signs on the matrix $\alpha $ then we will obtain
precisely the Markov transition matrix $A_{\mathcal{K}(f_{\mu })}$. In fact,
consider the $(n-1)\times (n-1)$ matrix 
\begin{equation*}
\beta :=\left[ 
\begin{array}{cc}
1_{n_{L}} & 0 \\ 
0 & -1_{n_{R}}
\end{array}
\right]
\end{equation*}
where $1_{n_{L}}$ and $1_{n_{R}}$ are the identity matrices of rank $n_{L}$
and $n_{R}$ respectively, with $n_{L}$ ($n_{R}$) being the number of
intervals $I_{i}$ of the Markov partition placed on the left (right) hand
side of the turning point of $f_{\mu }$. Therefore, we have 
\begin{equation*}
A_{\mathcal{K}(f_{\mu })}=\beta \alpha \text{.}
\end{equation*}
Now, consider the following matrix defined by 
\begin{equation*}
\gamma _{\mathcal{K}(f_{\mu })}:=(\gamma _{ij})\text{ \ \ with \ \ }\left\{ 
\begin{array}{l}
\gamma _{ii}=\varepsilon _{i}(c)\text{, \ \ }i=1,\ldots ,n \\ 
\\ 
\gamma _{in}=-\varepsilon _{i}(c)\text{, \ \ }i=1,\ldots ,n \\ 
\\ 
\gamma _{ij}=0\text{, \ \ otherwise.}
\end{array}
\right.
\end{equation*}
The matrix $\gamma _{\mathcal{K}(f_{\mu })}$ makes the diagram 
\begin{equation*}
\begin{array}{ccccc}
&  & \eta &  &  \\ 
& C_{0} & \longrightarrow & C_{1} &  \\ 
\gamma _{\mathcal{K}(f_{\mu })} & \downarrow &  & \downarrow & \beta \\ 
& C_{0} & \longrightarrow & C_{1} &  \\ 
&  & \eta &  & 
\end{array}
\end{equation*}
commutative. Finally, set $\theta _{\mathcal{K}(f_{\mu })}:=\gamma _{%
\mathcal{K}(f_{\mu })}\omega $. Then, the following diagram 
\begin{equation*}
\begin{array}{ccccc}
&  & \eta &  &  \\ 
& C_{0} & \longrightarrow & C_{1} &  \\ 
\theta _{\mathcal{K}(f_{\mu })} & \downarrow &  & \downarrow & A_{\mathcal{K}%
(f_{\mu })} \\ 
& C_{0} & \longrightarrow & C_{1} &  \\ 
&  & \eta &  & 
\end{array}
\end{equation*}
is also commutative. Now, notice that the transpose of $\eta $ has the
following factorization 
\begin{equation*}
\eta ^{T}=YiX\text{,}
\end{equation*}
where $Y$ is an invertible (over $\mathbb{Z}$) $n\times $ $n$ integer matrix
given by 
\begin{equation*}
Y:=\left( 
\begin{array}{cccccc}
1 & 0 & \cdots &  &  & 0 \\ 
0 & 1 & 0 & \cdots &  & 0 \\ 
\vdots & 0 & \ddots & \ddots &  & \vdots \\ 
& \vdots &  &  & 0 &  \\ 
0 & 0 & \cdots & 0 & 1 & 0 \\ 
-1 & -1 & \cdots &  & -1 & 1
\end{array}
\right) \text{,}
\end{equation*}
$i$ is the inclusion $C_{1}\hookrightarrow C_{0}$ given by 
\begin{equation*}
i:=\left( 
\begin{array}{cccc}
1 & 0 &  & 0 \\ 
0 & \ddots & \ddots & \vdots \\ 
\vdots & \ddots &  & 0 \\ 
&  &  & 1 \\ 
0 & \cdots &  & 0
\end{array}
\right)
\end{equation*}
and $X$ is an invertible (over $\mathbb{Z}$) $(n-1)\times (n-1)$ integer
matrix obtained from the $(n-1)\times n$\ matrix $\eta ^{T}$ by removing the 
$n$-th row of $\eta ^{T}$. Thus, from the commutative diagram 
\begin{equation*}
\begin{array}{ccccc}
&  & \eta ^{T} &  &  \\ 
& C_{1} & \longrightarrow & C_{0} &  \\ 
A_{\mathcal{K}(f_{\mu })}^{T} & \downarrow &  & \downarrow & \theta _{%
\mathcal{K}(f_{\mu })}^{T} \\ 
& C_{1} & \longrightarrow & C_{0} &  \\ 
&  & \eta ^{T} &  & 
\end{array}
\end{equation*}
we will have the following commutative diagram with short exact rows 
\begin{equation*}
\begin{array}{lllllllll}
&  &  & \text{ }i &  & \text{ }p &  &  &  \\ 
0 & \longrightarrow & C_{1} & \longrightarrow & C_{0} & \longrightarrow & 
C_{0}/C_{1} & \longrightarrow & 0 \\ 
&  & \downarrow A^{\prime } &  & \downarrow \theta ^{\prime } &  & 
\downarrow 0 &  &  \\ 
0 & \longrightarrow & C_{1} & \longrightarrow & C_{0} & \longrightarrow & 
C_{0}/C_{1} & \longrightarrow & 0 \\ 
&  &  & \text{ }i &  & \text{ }p &  &  & 
\end{array}
\end{equation*}
where the map $p$ is represented by the $1\times n$ matrix $\left[ 
\begin{array}{cccc}
0 & \ldots & 0 & 1
\end{array}
\right] $ and 
\begin{equation*}
A^{\prime }=XA_{\mathcal{K}(f_{\mu })}^{T}X^{-1}\text{ \ \ and \ \ }\theta
^{\prime }=Y^{-1}\theta _{\mathcal{K}(f_{\mu })}^{T}Y
\end{equation*}
i.e., $A^{\prime }$ is similar to $A_{\mathcal{K}(f_{\mu })}^{T}$\ over $%
\mathbb{Z}$ and $\theta ^{\prime }$ is similar to $\theta _{\mathcal{K}%
(f_{\mu })}^{T}$\ over $\mathbb{Z}$. Hence, for example by [\textbf{10}], we
obtain respectively 
\begin{eqnarray*}
\mathbb{Z}^{n-1}/(1-A^{\prime })\mathbb{Z}^{n-1} &\cong &\mathbb{Z}%
^{n-1}/(1-A_{\mathcal{K}(f_{\mu })})\mathbb{Z}^{n-1}\text{ and} \\
\mathbb{Z}^{n}/(1-\theta ^{\prime })\mathbb{Z}^{n} &\cong &\mathbb{Z}%
^{n}/(1-\theta _{\mathcal{K}(f_{\mu })})\mathbb{Z}^{n}\text{.}
\end{eqnarray*}
Now, from the last diagram we have, for example by [\textbf{9}], 
\begin{equation*}
\theta ^{\prime }=\left[ 
\begin{array}{cc}
A^{\prime } & \ast \\ 
0 & 0
\end{array}
\right] \text{.}
\end{equation*}
Therefore, 
\begin{equation*}
\mathbb{Z}^{n-1}/(1-A^{\prime })\mathbb{Z}^{n-1}\cong \mathbb{Z}%
^{n}/(1-\theta ^{\prime })\mathbb{Z}^{n}
\end{equation*}
and 
\begin{equation*}
\mathbb{Z}^{n-1}/(1-A_{\mathcal{K}(f_{\mu })})\mathbb{Z}^{n-1}\cong \mathbb{Z%
}^{n}/(1-\theta _{\mathcal{K}(f_{\mu })})\mathbb{Z}^{n}\text{.}
\end{equation*}
Next, we will compute $\mathbb{Z}^{n}/(1-\theta _{\mathcal{K}(f_{\mu })})%
\mathbb{Z}^{n}$. From the previous discussions and notations, the $n\times n$
matrix $\theta _{\mathcal{K}(f_{\mu })}$ is explicitly given by 
\begin{equation*}
\theta _{\mathcal{K}(f_{\mu })}:=\left( 
\begin{array}{ccccc}
-\varepsilon _{1}(c) & \varepsilon _{1}(c) & 0 & \cdots & 0 \\ 
\vdots & 0 & \ddots & \ddots & \vdots \\ 
& \vdots & \ddots &  & 0 \\ 
-\varepsilon _{n-1}(c) &  &  &  & \varepsilon _{n-1}(c) \\ 
0 & 0 & \cdots &  & 0
\end{array}
\right) \text{.}
\end{equation*}
Notice that this matrix $\theta _{\mathcal{K}(f_{\mu })}$ completely
describes the dynamics of $f_{\mu }$. Finally, using row and column
elementary operations over $\mathbb{Z}$, we can find invertible (over $%
\mathbb{Z}$) matrices $U_{1}$ and $U_{2}$ with integer entries such that 
\begin{equation*}
1-\theta _{\mathcal{K}(f_{\mu })}=U_{1}\left( 
\begin{array}{ccccc}
1+\sum\limits_{l=1}^{n-1}\prod\limits_{i=1}^{l}\varepsilon _{i}(c) &  &  & 
&  \\ 
& 1 &  &  &  \\ 
&  & \ddots &  &  \\ 
&  &  &  &  \\ 
&  &  &  & 1
\end{array}
\right) U_{2}\text{.}
\end{equation*}
Thus, we obtain 
\begin{equation*}
K_{0}(\mathcal{O}_{A_{\mathcal{K}(f_{\mu })}})\cong \mathbb{Z}^{n-1}/(1-A_{%
\mathcal{K}(f_{\mu })}^{T})\mathbb{Z}^{n-1}\cong \mathbb{Z}_{a}\text{,}
\end{equation*}
where $a=\left| 1+\sum\limits_{l=1}^{n-1}\prod\limits_{i=1}^{l}\varepsilon
_{i}(c)\right| $ and $n\in \mathbb{N}\backslash \left\{ 1\right\} .$
\end{proof}
\end{theorem}

\begin{example}
Set 
\begin{equation*}
\mathcal{K}(f_{\mu })=(RLLRRC)^{\infty },
\end{equation*}
where $R=-1,$ $L=+1,$ $C=0$. Thus, we can construct the $5\times 5$ Markov
transition matrix $A_{\mathcal{K}(f_{\mu })}$ and the matrices $\theta _{%
\mathcal{K}(f_{\mu })}$, $\omega $, $\varphi $, and $\pi $. 
\begin{equation*}
A_{\mathcal{K}(f_{\mu })}=\left( 
\begin{array}{ccccc}
0 & 1 & 1 & 0 & 0 \\ 
0 & 0 & 0 & 1 & 1 \\ 
0 & 0 & 0 & 0 & 1 \\ 
0 & 0 & 1 & 1 & 0 \\ 
1 & 1 & 0 & 0 & 0
\end{array}
\right) ,\text{ \ \ }\theta _{\mathcal{K}(f_{\mu })}=\left( 
\begin{array}{cccccc}
1 & -1 & 0 & 0 & 0 & 0 \\ 
-1 & 0 & 1 & 0 & 0 & 0 \\ 
-1 & 0 & 0 & 1 & 0 & 0 \\ 
1 & 0 & 0 & 0 & -1 & 0 \\ 
1 & 0 & 0 & 0 & 0 & -1 \\ 
0 & 0 & 0 & 0 & 0 & 0
\end{array}
\right) \text{,}
\end{equation*}
\begin{equation*}
\omega =\left( 
\begin{array}{cccccc}
0 & 1 & 0 & 0 & 0 & 0 \\ 
0 & 0 & 1 & 0 & 0 & 0 \\ 
0 & 0 & 0 & 1 & 0 & 0 \\ 
0 & 0 & 0 & 0 & 1 & 0 \\ 
0 & 0 & 0 & 0 & 0 & 1 \\ 
1 & 0 & 0 & 0 & 0 & 0
\end{array}
\right) \text{, \ \ }\varphi =\left( 
\begin{array}{cccccc}
-1 & 1 &  &  &  &  \\ 
& -1 & 1 &  &  &  \\ 
&  & -1 & 1 &  &  \\ 
&  &  & -1 & 1 &  \\ 
&  &  &  & -1 & 1
\end{array}
\right) \text{,}
\end{equation*}
\begin{equation*}
\pi =\left( 
\begin{array}{cccccc}
0 & 1 & 0 & 0 & 0 & 0 \\ 
0 & 0 & 1 & 0 & 0 & 0 \\ 
0 & 0 & 0 & 0 & 0 & 1 \\ 
0 & 0 & 0 & 1 & 0 & 0 \\ 
0 & 0 & 0 & 0 & 1 & 0 \\ 
1 & 0 & 0 & 0 & 0 & 0
\end{array}
\right) \text{.}
\end{equation*}
We have 
\begin{equation*}
K_{0}(\mathcal{O}_{A_{\mathcal{K}(f_{\mu })}})\cong \mathbb{Z}_{2}\text{ \ \
and \ \ }K_{1}(\mathcal{O}_{A_{\mathcal{K}(f_{\mu })}})\cong \{0\}.
\end{equation*}
\end{example}

\begin{remark}
In the statement of Theorem 1 the case $a=0$ may occur. This happens when we
have a star product factorizable kneading sequence [\textbf{4}]. In this
case the correspondent Markov transition matrix is reducible.
\end{remark}

\begin{remark}
In [\textbf{6}] the authors have constructed a class of $C^{\ast }$-algebras
from the $\beta $-expansions of real numbers. In fact, by considering a
semiconjugacy from the real quadratic map to the tent map [\textbf{14}], we
can also obtain Theorem 1 using [\textbf{6}] and the $\lambda $-expansions
of real numbers introduced in [\textbf{4}].
\end{remark}

\begin{remark}
In [\textbf{12}] (see also [\textbf{11}])\ and [\textbf{13}] the BF-groups
are explicitly calculated with respect to other kind of maps on the
interval. \bigskip
\end{remark}

\end{document}